\documentclass{amsproc}
\usepackage{amssymb,amscd}
\usepackage[cp1251]{inputenc}
\usepackage[T2A]{fontenc}
\usepackage[matrix,arrow]{xy}

\newtheorem{theorem}{Theorem}
\newtheorem{definition}{Definition}
\newtheorem{lemma}{Lemma}
\newtheorem{proposition}{Proposition}
\newtheorem{remark}{Remark}

\newcommand{\bideg}{\mathrm{bideg}}

\newcommand{\supp}{\mathrm{supp}}
\newcommand{\im}{\mathrm{Im}}

\begin{document}
\title[The Hodge filtration and integral representations]{The Hodge filtration on complements of complex coordinate subspace arrangements and integral representations of holomorphic functions}

\author{Yury~V. Eliyashev}
\address{Institute of Mathematics and Computer Science,
Siberian Federal University, Svobodny 79, 660041 Krasnoyarsk,
Russia} \email{eliyashev@mail.ru}

\begin{abstract}
We compute the Hodge filtration on cohomology groups of complements
of complex coordinate subspace arrangements. By means of this result
we construct integral representations of holomorphic functions such
that kernels of these representations have singularities on complex
coordinate subspace arrangements.
\end{abstract}

\maketitle

\section*{Introduction}
A study of topology of coordinate subspace arrangements appears in
different areas of mathematics: in toric topology and combinatorial
topology \cite{BP,BP2}, in the theory of toric varieties, where
complements to  coordinate subspace arrangements play the role of
homogeneous coordinate spaces \cite{CX,CX2}, in the theory of
integral representations of holomorphic functions in several complex
variables, where coordinate subspace arrangements play the role of
singular sets of integral representations kernels \cite{AJ,Sh}.

The universal combinatorial method for the computation of cohomology
groups of complements to \emph{arbitrary} subspace arrangements was
developed in the book of Goresky and Macpherson \cite{GM} (see also
\cite{VS}), but this method often leads to cumbersome computations.
In the study of toric topology, in particular, in works of
Buchstaber and Panov \cite{BP, BP2}, the method for the computation
of the cohomology of complements to \emph{coordinate} subspace
arrangements was developed, this method is simpler than the
universal method and allows to get some additional topological
information.

The main purpose of this article is to compute the Hodge filtration
on the cohomology rings of complements to complex coordinate
subspace arrangements. We will show that
 the Hodge filtration is described by means of a special bigrading on the cohomology rings of complements to complex coordinate subspace arrangements, which was introduced in \cite{BP,BP2}, this bigrading was obtained originally from the combinatorial and topological ideas. We use these results to construct the integral representations of holomorphic functions such that the kernels of these representations have singularities on coordinate subspace arrangements.

The first section of this paper consists of different facts about
topology of complements to complex coordinate subspace arrangements,
in the text of this section we follow \cite{BP}, \cite{BP2}.  Let
$Z$ be a complex coordinate subspace arrangement in $\mathbb{C}^n.$
In \cite{BP}, \cite{BP2}, from the topological reasons, the
differential bigraded  algebra $R$ was introduced ($R$ is determined
by combinatorics of $Z$) such that the ring of  cohomology
$H^*(\mathbb{C}^n \setminus Z)$ is isomorphic to the ring of
cohomology $H^*(R).$ Denote by $H^{p,q}(R)$ the bigraded cohomology
of the algebra $R,$ then $$H^s(\mathbb{C}^n \setminus Z)\simeq
\bigoplus_{p+q=s} H^{p,q}(R).$$ Thus, we get a bigrading on the
cohomology ring $H^*(\mathbb{C}^n \setminus Z)$.

In the second section we recall some facts and  concept from
differential topology and complex analysis. These facts we use in
the last two sections.

In the third section the main theorem of this paper is proved. We
will show that the bigrading on the cohomology of $R$ and,
consequently, the bigrading on the  cohomology
$H^*(\mathbb{C}^n\setminus Z)$ appear naturally from the complex
structure on the manifold $\mathbb{C}^n\setminus Z.$ In particular,
denote by $F^k H^s(\mathbb{C}^n\setminus Z,\mathbb{C})$ a $k$-th
term of the Hidge filtration on $H^s(\mathbb{C}^n\setminus
Z,\mathbb{C}).$ Then there is the following theorem.
\begin{theorem}
$$F^k H^s(\mathbb{C}^n\setminus Z,\mathbb{C})=\bigoplus_{p\geq k} H^{p,s-p}(R,\mathbb{C}).$$
\end{theorem}

In the last section we construct integral representations of
holomorphic functions such that kernels of these representations
have singularities on coordinate subspace arrangements.

\section{General facts on topology of coordinate subspace arrangements}
In this section different facts about topology of complements to
coordinate subspaces arrangements are gathered. All statements of
this section are taken from \cite{BP}.

Let $\mathcal{K}$ be an arbitrary simplicial complex on the vertex
set $[n]=\{1,\dots,n\}.$ Define a coordinate planes arrangement
$$Z_{\mathcal{K}}:=\bigcup_{\sigma\not\in\mathcal{K}}L_\sigma,$$ where $\sigma=\{i_1,\dots,i_m\}\subseteq[n]$ is a subset in $[n]$ such that $\sigma$ does not define a simplex in $\mathcal{K}$ and $$L_\sigma=\{z\in\mathbb{C}^n: z_{i_1}=\dots=z_{i_m}=0\}.$$
Any arrangement of complex coordinate subspaces in $\mathbb{C}^n$ of
codimension greater than $1$ can be defined in this way.

Consider a cover
$\mathcal{U}_\mathcal{K}=\{\mathcal{U}_\sigma\}_{\sigma\in
\mathcal{K}}$ of $\mathbb{C}^n\setminus Z_\mathcal{K},$ where
$$\mathcal{U}_\sigma=\mathbb{C}^n\setminus\bigcup_{i\not\in \sigma}\{z_i=0\}.$$

By $D^2_\sigma\times S^1_\gamma$ denote the following chain
$$D^2_\sigma\times S^1_\gamma=\{|z_i|\leq 1:i\in \sigma; |z_j|=1: j\in\gamma,z_k=1: k\not\in\gamma\cup\sigma\},$$
where $\sigma,\gamma\subseteq[n]$ and $\sigma\cap\gamma=\emptyset.$
We define the form
 \begin{equation}\label{eq.dz}\frac{d z_I}{z_I}=\frac{d z_{i_1}}{z_{i_1}}\wedge\dots \wedge \frac{d z_{i_{k}}}{z_{i_{k}} },\end{equation}
 where $I\subseteq[n], |I|=k,$ $I=\{i_1,\dots,i_{k}\},$ and $i_1<\dots<i_{k}.$

 The orientation of the chain $D^2_\sigma\times S^1_\gamma$ is such that the restriction of the form
$$\frac{1}{(\sqrt{-1})^{|\gamma|}}\frac{d z_\gamma}{z_\gamma} \wedge \bigwedge_{j\in \sigma}(\sqrt{-1} d z_j\wedge d \overline{z}_j)$$
 on $D^2_\sigma\times S^1_\gamma$ is positive. Then the boundary of this chain equals
$$\partial D^2_\sigma\times S^1_\gamma=\Sigma_{i\in\sigma}(-1)^{(i,\gamma)}D^2_{\sigma\setminus i}\times S^1_{\gamma\cup i},$$
where $(i,\gamma)$ is the position of $i$ in the naturally ordered
set $\gamma\cup i.$

\begin{definition}\label{def.MA-comlex}The topological space
$$\mathcal{Z}_\mathcal{K}=\bigcup_{\sigma\in \mathcal{K}} D^2_\sigma\times S^1_{[n]\setminus \sigma}$$
is called a moment-angle complex.
\end{definition}

\begin{theorem}[\cite{BP}]\label{th.ret}
There is a deformation retraction from $\mathbb{C}^n\setminus
Z_\mathcal{K}$ to $\mathcal{Z}_\mathcal{K}.$
\end{theorem}

\begin{definition}\label{def.SR-ring}
 A Stanley–Reisner ring of a simplicial complex  $\mathcal{K}$ on the vertex set $[n]$ is a ring
  $$\mathbb{Z}[\mathcal{K}]=\mathbb{Z}[v_1,\dots,v_n]/ \mathcal{I}_\mathcal{K},$$ where $\mathcal{I}_\mathcal{K}$ is a homogeneous ideal generated by the monomials $v_\sigma=\prod_{i\in\sigma}v_i$ such that $\sigma\not\in\mathcal{K}:$ $$\mathcal{I}_\mathcal{K}=(v_{i_1}\cdot{\dots}\cdot v_{i_m}: \{i_1,\dots,i_m\}\not\in\mathcal{K}).$$
\end{definition}

Consider a differential bigraded algebra
$(R(\mathcal{K}),\delta_R):$
$$R_\mathcal{K}:=\Lambda[u_1,\dots,u_n]\otimes
\mathbb{Z}[\mathcal{K}]/\mathcal{J},$$ where
$\Lambda[u_1,\dots,u_n]$ is an exterior algebra, $\mathcal{J}$ is
the ideal generated by monomials $v_i^2,u_i\otimes v_i,
i=1,\dots,n.$ Bidegrees of generators $v_i, u_i$  of this algebra
are equal to
$$\bideg\; v_i=(1,1), \bideg\; u_i=(1,0).$$
The differential $\delta_R$ is defined on the generators as follows
$$\delta_R u_i=v_i, \delta_R v_i =0.$$
\begin{remark}
In \cite{BP} a \textbf{different bigrading} on the algebra
$R_\mathcal{K}$ was used, but our bigrading is equivalent to the
bigrading from \cite{BP}.
\end{remark}

We denote by $R^{p,q}_\mathcal{K}$ the homogeneous component of the
algebra $R_\mathcal{K}$ of the bidegree $(p,q).$  The differential
$\delta_R$ is compatible with the bigrading, i.e.,
$\delta_R(R^{p,q}_\mathcal{K})\subseteq R^{p,q+1}_\mathcal{K}.$
Consider the complex
$$\dots\stackrel{\delta_R}{\longrightarrow}R^{p,q-1}_\mathcal{K}\stackrel{\delta_R}{\longrightarrow}R^{p,q}_\mathcal{K}\stackrel{\delta_R}{\longrightarrow}R^{p,q+1}_\mathcal{K}\stackrel{\delta_R}{\longrightarrow}\dots,$$
denote by $H^{p,q}(R_\mathcal{K})$ a cohomology group of this
complex.  It is clear that the cohomology of $R_\mathcal{K}$ are
isomorphic to
$$H^s(R_\mathcal{K})=\bigoplus_{p+q=s}H^{p,q}(R_\mathcal{K}).$$

\begin{theorem}[\cite{BP}]\label{th.BP1}
The cohomology ring $H^*(\mathbb{C}^n\setminus Z_{\mathcal{K}})$ is
isomorphic to the ring $H^*(R_\mathcal{K}).$
\end{theorem}

\begin{remark}
The relation between Theorem \ref{th.BP1} and the results of Goresky
and Macpherson \cite{GM} on cohomology of subspace arrangements is
described in \cite[Ch.~8]{BP}.
\end{remark}

Now we describe the explicit construction of the isomorphism of
Theorem \ref{th.BP1}. First, we construct a cell decomposition of
$\mathcal{Z}_\mathcal{K}.$ Define a cell
$$E_{\sigma\gamma}=\{|z_i|<1:i\in \sigma; |z_j|=1,z_j\neq1: j\in\gamma;z_k=1: k\not\in\gamma\cup\sigma\},$$
where $\sigma,\gamma\subseteq[n]$ and $\sigma\cap\gamma=\emptyset.$
The closure of this cell equals
$\overline{E_{\sigma\gamma}}=D^2_\sigma\times S^1_\gamma.$ The
orientation of $E_{\sigma\gamma}$  is defined by the orientation of
$D^2_\sigma\times S^1_\gamma.$ We obtain the cell decomposition
$$\mathcal{Z}_\mathcal{K} = \bigcup_{\sigma\in \mathcal{K}, \gamma\subseteq [n]\setminus \sigma}E_{\sigma\gamma}.$$

Let $C_*(\mathcal{Z}_\mathcal{K})$ be the group of cell chains of
this cell decomposition, then denote by
$C^*(\mathcal{Z}_\mathcal{K})$ the group of cell cochains. Let
$E_{\sigma\gamma}'$ be  a cocell dual to the cell
$E_{\sigma\gamma},$ i.e., $E_{\sigma\gamma}'$ is a linear functional
from  $C_*(\mathcal{Z}_\mathcal{K})$ such that $\langle
E_{\sigma\gamma}',E_{\sigma'\gamma'}\rangle=\delta^{\sigma\gamma}_{\sigma'\gamma'}$
(the Kronecker delta).

Denote $u_I v_J:=u_{i_1}\dots u_{i_q}\otimes v_{j_1}\dots v_{j_p},$
where $I=\{i_1,\dots,i_q\},$
 $i_1<\dots<i_q,$ $J=\{j_1,\dots,j_p\},$ and $I\cap J=\emptyset,I,J\subseteq[n],$ (we suppose that $u_\emptyset v_\emptyset=1$).
\begin{proposition}[\cite{BP}]\label{th.BP2}
The linear map $\phi:R_\mathcal{K}\rightarrow
C^*(\mathcal{Z}_\mathcal{K}),$ $\phi(v_\sigma u_\gamma) =
E_{\sigma\gamma}'$ is an isomorphism of differential bigraded
modules. In particular, there is an additive isomorphism
$H^*(R_\mathcal{K})\stackrel{\phi}{\simeq}
H^*(\mathcal{Z}_\mathcal{K}).$
\end{proposition}
From the structure of the cell decomposition of
$\mathcal{Z}_\mathcal{K}$ and Theorem \ref{th.ret} we obtain that
every cycle $\Gamma\in H_s(\mathbb{C}^n\setminus Z_\mathcal{K})$ has
a representative of the form
\begin{equation}\label{eq.cycle_decom}
\Gamma=\sum_{p+q=s}\Gamma_{p,q},
\end{equation}
where $\Gamma_{p,q}$ is a cycle of the form
\begin{equation}\label{eq.cycle_type}\Gamma_{p,q}=\sum_{\substack{|\sigma|=q\\|\gamma|=p-q}}C_{\sigma \gamma}\cdot D^2_\sigma\times S^1_\gamma.\end{equation}
 A group generated by all cycles of the form  (\ref{eq.cycle_type}) is denoted by  $H_{p,q}(\mathbb{C}^n\setminus Z_{\mathcal{K}}).$ Obviously we have  $$H_{s}(\mathbb{C}^n\setminus Z_{\mathcal{K}})=\bigoplus_{p+q=s}H_{p,q}(\mathbb{C}^n\setminus Z_{\mathcal{K}}).$$

It follows from Proposition \ref{th.BP2} that
$\langle\Gamma_{p,q},\phi(\omega^{p',q'})\rangle=0$ for any
$\Gamma_{p,q}\in H_{p,q}(\mathbb{C}^n\setminus~Z_{\mathcal{K}})$ and
$\omega^{p',q'}\in H^{p',q'}(R_\mathcal{K}),$ $p'\neq p$ and $q'\neq
q.$ Hence, the pairing between $\Gamma_{p,q}$ and
$\phi(\omega^{p',q'})$ can be nonzero only if $p'=p,q'=q.$ Therefore
the pairing  between the vector spaces
$H_{p,q}(\mathbb{C}^n\setminus Z_{\mathcal{K}},\mathbb{R})$ and
$\phi(H^{p',q'}(R_\mathcal{K}\otimes\mathbb{R}))$ is nondegenerate
if $p=p',$ $q=q'$ and equals to zero otherwise.

\section{Cech cohomology, filtrations and cochains}
In this section we recall some facts from differential topology and
complex analysis, we mainly use a material from the books \cite{BT},
\cite{GH}. Let $X$ be a complex manifold and
$\mathcal{U}=\{\mathcal{U}_\alpha\}_{\alpha\in \mathcal{A}}$ is an
open, countable, locally finite cover of this manifold. Now we
introduce the following notation for sheafs on $X$: $\mathcal{E}^s$
denotes the sheaf of {$C^\infty$-differential} forms of degree $s$,
$\mathcal{E}^{p,q}$ denotes the sheaf of $C^\infty$-differential
forms of bedegree $(p,q)$, $\Omega^p$ denotes the sheaf of
holomorphic differential forms of degree $p.$
\begin{definition}
The decreasing filtration $$F^k
\mathcal{E}^{\bullet}=\bigoplus_{p\geq
k}\mathcal{E}^{p,\bullet-p},$$ on the de Rham complex
$(\mathcal{E}^\bullet,d)$ is called
 the Hodge filtration.
\end{definition}

The Hodge filtration induces a filtration $F^k H^s(X,\mathbb{C})$ on
a de Rham cohomology, i.e.,
$$F^k H^s(X,\mathbb{\mathbb{C}})=\im(H^s(F^k \mathcal{E}^\bullet(X),d)\rightarrow H^s(\mathcal{E}^\bullet(X),d)),$$
 where $H^s(\mathcal{E}^\bullet(X),d)$ is the cohomology of the de Rham complex and $H^s(F^k \mathcal{E}^\bullet(X),d)$ is the cohomology of $k$-th term of the Hodge filtration. In other words, if $\omega$ lies in $F^k H^s(X,\mathbb{C})$ then there is a form $\widetilde{\omega}, [\widetilde{\omega}]=\omega$ such that $$\widetilde{\omega}=\sum_{\substack{p\geq k}}\widetilde{\omega}^{p,s-p},$$
where $\widetilde{\omega}^{p,q}\in \mathcal{E}^{p,q}(X).$

Let $C^t(\mathcal{E}^s,\mathcal{U})$ be the Cech-de Rham double
complex for the cover $\mathcal{U}$:
$C^t(\mathcal{E}^s,\mathcal{U})$ with a Cech coboundary operator
$\delta: C^t(\mathcal{E}^s,\mathcal{U})\rightarrow
C^{t+1}(\mathcal{E}^s,\mathcal{U})$ and a de Rham differential $d:
C^t(\mathcal{E}^s,\mathcal{U})\rightarrow
C^{t}(\mathcal{E}^{s+1},\mathcal{U})$ on this complex, i.e.,
$$(\delta \omega)_{i_0,\dots,i_{t+1}}=(-1)^s\sum_{j=0}^{t+1}(-1)^j \omega_{i_0,\dots,\widehat{i_j},\dots,i_{t+1}}|_{\mathcal{U}_{i_0}\cap\dots\cap\mathcal{U}_{i_{t+1}}},$$
$$(d\omega)_{i_0,\dots,i_{t}}=d(\omega)_{i_0,\dots,i_{t}}.$$

The associated single complex is defined by
$$K^r(\mathcal{U},\mathcal{E}^\bullet)=\bigoplus_{s+t=r}C^t(\mathcal{E}^s,\mathcal{U})$$
the operator $D=\delta+d$ is the differential of this complex.
Notice that our definition of Cech coboundary $\delta$ is different
from the standard one by the factor $(-1)^s$, with this choice of
sign we get $D^2=0,$ hence
$(K^\bullet(\mathcal{U},\mathcal{E}^\bullet),D)$ is a complex. There
is a natural inclusion of the de Rham complex
$\varepsilon:\mathcal{E}^{\bullet}(X)\rightarrow
C^0(\mathcal{E}^\bullet,\mathcal{U}),$
$\varepsilon(\omega)_{j_0}=\omega|_{\mathcal{U}_{j_0}},$ also we
denote the induced map from $\mathcal{E}^{\bullet}(X)$ to
$K^\bullet(\mathcal{U},\mathcal{E}^\bullet)$ by $\varepsilon.$

\begin{theorem}{\cite{BT}}\label{th.BT1}
The inclusion $\varepsilon:\mathcal{E}^{\bullet}(X)\rightarrow
K^\bullet(\mathcal{U})$ is a quasi-isomorphism of complexes, i.e.,
$H^s(X,\mathbb{C})\stackrel{\varepsilon}{\simeq}
H^s(K^\bullet(\mathcal{U},\mathcal{E}^\bullet),D).$
\end{theorem}

 The Hodge filtration $F^k K^\bullet(\mathcal{U},\mathcal{E}^\bullet)$  is defined naturally on $(K^\bullet(\mathcal{U},\mathcal{E}^\bullet),D)$. This filtration induces a filtration on cohomology $F^k H^s(K^\bullet(\mathcal{U},\mathcal{E}^\bullet),D).$
There is an isomorphism $F^k
H^s(X,\mathbb{C})\stackrel{\varepsilon}{\simeq} F^k
H^s(K^\bullet(\mathcal{U},\mathcal{E}^\bullet),D).$

Consider a subcomplex $K^r(\mathcal{U},\Omega^\bullet)$ of the
complex $K^r(\mathcal{U},\mathcal{E}^\bullet)$
$$K^r(\mathcal{U},\Omega^\bullet)=\bigoplus_{s+t=r}C^t(\Omega^s,\mathcal{U}),$$
and an  inclusion map $\tau: K^r(\mathcal{U},\Omega^\bullet)
\rightarrow K^r(\mathcal{U},\mathcal{E}^\bullet).$ It is easy to get
the following statement.
\begin{theorem}\label{th.BT2}
Suppose $\mathcal{U}$ is a $\overline{\partial}$-acyclic cover of
$X$ then the inclusion  $\tau$ is  a quasi-isomorphism of the
complexes $K^\bullet(\mathcal{U},\Omega^\bullet)$ and
$K^\bullet(\mathcal{U},\mathcal{E}^\bullet).$
\end{theorem}

Let $F^k \Omega^p$ be a stupid filtration on the de Rham complex of
holomorphis forms $(\Omega^\bullet,d),$ i.e.,
\begin{equation*}
F^k \Omega^p=\left\{
           \begin{array}{ll}
             \Omega^p & \mbox{for\;} p\geq k.\\
             0 & \mbox{for\;} p<k. \\
           \end{array}
         \right. \end{equation*}

The stupid filtration induces filtration on cohomology $F^k
H^s(K^\bullet(\mathcal{U},\Omega^\bullet),D).$ Suppose $\mathcal{U}$
is a $\overline{\partial}$-acyclic cover of $X$ then $F^k
H^s(K^\bullet(\mathcal{U},\Omega^\bullet),D)\simeq F^k
H^s(X,\mathbb{C}).$

From now until the end of this section we will follow the paper
\cite{Gleason}.
\begin{definition}
A $\mathcal{U}$-chain of degree $t$ and of dimension $s$ on the
manifold $X$ is an alternating function $\Gamma$ from the set of
indexes $\mathcal{A}^{t+1}$ to the group of singular chains in $X$
of dimension $s$ such that $\Gamma$ is nonzero on a finite number of
points from $\mathcal{A}^{t+1}$ and
$$\supp(\Gamma_{i_0,\dots,i_t})\subset
\mathcal{U}_{i_0}\cap\dots\cap\mathcal{U}_{i_t},$$  for every
$(i_0,\dots,i_{t})\in\mathcal{A}^{t+1},$ where
$\supp(\Gamma_{i_0,\dots,i_t})$ is the support of the chain
$\Gamma_{i_0,\dots,i_t}.$
\end{definition}

Let  $C_{t,s}(\mathcal{U})$ be an additive group of
$\mathcal{U}$-chains of degree $t$ and of dimension $s$ on the
manifold $X.$ Define maps $\delta':C_{t,s}(\mathcal{U})\rightarrow
C_{t-1,s}(\mathcal{U})$
$$(\delta'\Gamma)_{i_0,\dots,i_{t-1}}= (-1)^s\sum_{i\in\mathcal{A}}\Gamma_{i,i_0,\dots,i_{t-1}},$$
and $\partial:C_{t,s}(\mathcal{U})\rightarrow
C_{t,s-1}(\mathcal{U})$
$$(\partial\Gamma)_{i_0,\dots,i_{t}}=\partial(\Gamma)_{i_0,\dots,i_{t}},$$
i.e., the operator $\partial$ is a boundary operator on each chain
$\Gamma_{i_0,\dots,i_{t}}.$ The groups $C_{t,s}(\mathcal{U}),t,s\geq
0$ together with the differentials $\delta',\partial$ form a double
complex. Define a map $\varepsilon':C_{0,s}(\mathcal{U})\rightarrow
C_{s}(X)$ in the following way
$$\varepsilon'(\Gamma)=\sum_{i\in\mathcal{A}}\Gamma_i.$$

Now we will construct a pairing between elements of
$C_{t,s}(\mathcal{U})$ and $C^t(\mathcal{E}^s,\mathcal{U}).$ Suppose
$\Gamma \in C_{t,s}(\mathcal{U})$ and $\omega\in
C^t(\mathcal{E}^s,\mathcal{U}),$ then
$$\langle\omega,\Gamma\rangle= \frac{1}{(t+1)!}\sum_{(i_0,\dots,i_t)\in \mathcal{A}^{t+1}}\int_{\Gamma_{i_0,\dots,i_t}}\omega_{i_0,\dots,i_t}.$$

There are the following relations for the pairing:
\begin{gather*}\label{eq.pairing}
\langle\omega^{t,s},\partial\Gamma_{t,s+1}\rangle=\langle d\omega^{t,s},\Gamma_{t,s+1}\rangle,\\
\langle\delta\omega^{t,s},\Gamma_{t+1,s}\rangle=\langle \omega^{t,s},\delta'\Gamma_{t+1,s}\rangle,\\
\int_{\varepsilon'(\Gamma_{0,s})}\omega^s=\langle\varepsilon\omega^s,\Gamma_{0,s}\rangle,
\end{gather*}
where $\omega^{t,s}\in C^t(\mathcal{E}^s,\mathcal{U})$, $\omega^s\in
\mathcal{E}^s(X),$ and $\Gamma_{t,s} \in C_{t,s}(\mathcal{U}).$

\begin{definition}
Let $\Gamma$ be a singular cycle of dimension $s$ on $X,$ then a
$\mathcal{U}$-resolvent of length $k$ of the cycle $\Gamma$ is a
collection of  $\mathcal{U}$-chains $\Gamma^i\in
C_{i,s-i}(\mathcal{U}),$ $0\leq~i\leq~k$ such that
$\Gamma=\varepsilon' \Gamma^0$ and $\partial \Gamma^i=-\delta'
\Gamma^{i+1}.$
\end{definition}

\begin{proposition}\label{prop.pairing}
Given an $s$-dimensional cycle $\Gamma,$ a closed differential form
$\omega$ of degree $s,$ a $\mathcal{U}$-resolvent
$\Gamma^0,\dots,\Gamma^k$ of the cycle $\Gamma$ and a cocycle
$\widetilde{\omega}\in K^s(\mathcal{U})$ such that
$\widetilde{\omega}=\sum_{i\leq k}{\widetilde{\omega}^{i,s-i}},$
$\widetilde{\omega}^{i,s-i}\in C^{i}(\mathcal{E}^{s-i},\mathcal{U})$
and the cocycle $\varepsilon \omega$ is cohomologous to
$\widetilde{\omega}$ in
$H^s(K^\bullet(\mathcal{U},\mathcal{E}^\bullet),D),$ then
$$\int_{\Gamma}\omega=\sum_{i\leq k}\langle\widetilde{\omega}^{i,s-i},\Gamma^i\rangle.$$
\end{proposition}
This proposition follows directly from the properties of the
pairing.

\section{The Hodge filtration of cohomology of complements to coordinate subspace arrangements}
In this section we compute the Hodge filtration of the cohomology
ring $H^*(\mathbb{C}^n\setminus~Z_{\mathcal{K}},\mathbb{C}).$ It
follows from Theorem \ref{th.BP1} and Proposition \ref{th.BP2} that
there is the isomorphism $H^*(\mathbb{C}^n\setminus
Z_{\mathcal{K}},\mathbb{C})\stackrel{\phi}{\simeq}H^*(R_{\mathcal{K}}\otimes\mathbb{C}).$
\setcounter{theorem}{0}
\begin{theorem}\label{th.main}Let $H^{p,q}(R_{\mathcal{K}}\otimes\mathbb{C})$  be the bigraded cohomology group of the complex $R^{p,q}_{\mathcal{K}}\otimes\mathbb{C},$ then there is an isomorphism
$$F^k H^s(\mathbb{C}^n\setminus Z_{\mathcal{K}},\mathbb{C})\stackrel{\phi}{\simeq}\bigoplus_{p\geq k}H^{p,s-p}(R_{\mathcal{K}}\otimes\mathbb{C}).$$
\end{theorem}
\setcounter{theorem}{5} {\it Proof.} First, we will prove the lemma.
\begin{lemma}\label{lemma.resolv} Let $$\Gamma_{p,q}=\sum_{\substack{|\sigma|=q\\|\gamma|=p-q}}C_{\sigma \gamma}\cdot D^2_\sigma\times S^1_\gamma$$ be a cycle in $\mathbb{C}^n\setminus Z_{\mathcal{K}}.$ Then there is a $\mathcal{U}_\mathcal{K}$-resolvent of the cycle $\Gamma_{p,q}$ of length $q$:  $$\Gamma^{0}_{p,q},\dots,\Gamma^{q}_{p,q},$$  where $\Gamma^{k}_{p,q}$ is a $\mathcal{U}_\mathcal{K}$-chain of dimension $q+p-k$ and of degree $k$ of the form
$$(\Gamma^{k}_{p,q})_{\alpha_0,\dots,\alpha_k}=\sum_{\substack{|\sigma|=q-k\\|\gamma|=p-q+k}}C_{\sigma \gamma, \alpha_0\dots\alpha_k}\cdot D^2_\sigma\times S^1_\gamma.$$
\end{lemma}
{\it Proof.} We will use the induction on the length $k$ of the
resolvent. We going to construct the resolvent of the special form
$$(\Gamma^{k}_{p,q})_{\sigma_k,\sigma_{k-1},\dots,\sigma_0}=\sum_{|\gamma|=p-q+k}C_{\sigma_k \gamma, \sigma_k\dots \sigma_0}\cdot D^2_{\sigma_k}\times S^1_\gamma,$$
for $|\sigma_j|=q-j,$  $\sigma_{j}\subset \sigma_t, j>t$
$j=0,\dots,k$ (in other words, $\{\sigma_j\}$ is a chain of subsets
$\sigma_k\subset\dots\subset\sigma_0,$ and $|\sigma_j|=q-j$); and
$(\Gamma^{k}_{p,q})_{\alpha_0,\dots,\alpha_k}=0$ for any other
indexes $\alpha_0,\dots,\alpha_k.$

The base of induction: define
$(\Gamma^0_{p,q})_{\sigma_0}=\sum_{\substack{|\gamma|=p-q}}C_{\sigma_0
\gamma}\cdot D^2_{\sigma_0}\times S^1_\gamma$ with $|\sigma_0|=q$
and $(\Gamma^0_{p,q})_{\alpha}=0$ for any other indexes $\alpha.$ We
get
$$\Gamma_{p,q}=\sum_{\substack{|\sigma_0|=q\\|\gamma|=p-q}}C_{\sigma_0 \gamma}\cdot D^2_{\sigma_0}\times S^1_\gamma=\sum_{\sigma\in\mathcal{K}}(\Gamma^{0}_{p,q})_{\sigma}=\varepsilon' \Gamma^{0}_{p,q},$$
therefore $\Gamma^{0}_{p,q}$ is the resolvent of length $0$.

Suppose that the resolvent $\Gamma^{0}_{p,q},\dots,\Gamma^{k}_{p,q}$
of length $k$ is already constructed. Recall that $(i,\gamma)$  is
the position of $i$ in the naturally ordered set $\gamma\cup i.$
Define $$(\Gamma^{k+1}_{p,q})_{\sigma_k\setminus i,\sigma_k
\dots\sigma_0} =(-1)^{p+q-k}\sum_{\substack{|\gamma|=p-q}}
(-1)^{(i,\gamma)}C_{\sigma_k \gamma, \sigma_k \dots\sigma_0}
D^2_{\sigma_k\setminus i}\times S^1_{\gamma\cup i},$$ for $i\in
\sigma_k,$ $|\sigma_j|=q-j,$ $\sigma_{j+1}\subset \sigma_{j},$ and
$$(\Gamma^{k+1}_{p,q})_{\alpha_0,\dots,\alpha_{k+1}}=0,$$ for any
other indexes $\alpha_0,\dots,\alpha_{k+1}.$ Let us show that
$\Gamma^{0}_{p,q},\dots,\Gamma^{k+1}_{p,q}$ is a resolvent of length
$k+1:$
\begin{multline*}-(\delta' \Gamma^{k+1}_{p,q})_{\sigma_k \dots\sigma_0}=(-1)^{p+q-k}\sum_{i\in \sigma_k}(\Gamma^{k+1}_{p,q})_{\sigma_k\setminus i,\sigma_k \dots\sigma_0}=\\=\sum_{\substack{i\in \sigma_k\\|\gamma|=p-q}} (-1)^{(i,\gamma)}C_{\sigma_k \gamma, \sigma_k \dots\sigma_0} D^2_{\sigma_k\setminus i}\times S^1_{\gamma\cup i}
=\sum_{|\gamma|=p-q}C_{\sigma_k \gamma, \sigma_k \dots\sigma_0}
\partial D^2_{\sigma_k}\times S^1_{\gamma}=\\=(\partial
\Gamma^{k}_{p,q})_{\sigma_k \dots\sigma_0}.\end{multline*}
 For any indexes $\alpha_0,\dots,\alpha_{k}$ different from $\sigma_k,\dots,\sigma_{m+1},\sigma_{m-1},\dots,\sigma_0,$ $0\leq m \leq k,$  directly from definition of $\Gamma^k_{p,q},\Gamma^{k+1}_{p,q}$, we get $$(\partial\Gamma^{k}_{p,q})_{\alpha_0,\dots,\alpha_{k}}=-(\delta'\Gamma^{k+1}_{p,q})_{\alpha_0,\dots,\alpha_{k}}=0.$$
Consider the last case $\sigma_k\setminus
i,\sigma_k,\dots,\sigma_{m+1},\sigma_{m-1},\dots,\sigma_0,$ for
$0\leq m \leq k.$ Since by the induction hypothesis
$\Gamma^{0}_{p,q}\dots\Gamma^{k}_{p,q}$  is a resolvent,
$-\delta'\Gamma^{k}_{p,q}=\partial \Gamma^{k-1}_{p,q},$ hence we
have $\delta' \partial \Gamma^{k}_{p,q}=0,$ and
\begin{multline*}(\delta' \partial\Gamma^{k}_{p,q})_{\sigma_k\dots\sigma_{m+1}\sigma_{m-1}\dots\sigma_0}=\\=(-1)^{p+q-k+1}
\sum_{\substack{\sigma_{m+1}\subset\sigma_m\subset\sigma_{m-1}\\|\sigma_m|=q-m}}\sum_{|\gamma|=p-q+k}\sum_{i\in
\sigma_k} (-1)^{(i,\gamma)} C_{\sigma_k \gamma,
\sigma_k\dots\sigma_{0}}\cdot D^2_{\sigma_k\setminus i}\times
S^1_{\gamma\cup i}=0.\end{multline*} Therefore, for a fixed $i\in
\sigma_k,$ we get
$$\sum_{\substack{\sigma_{m+1}\subset\sigma_m\subset\sigma_{m-1}\\|\sigma_m|=q-m}}\sum_{|\gamma|=p-q+k} (-1)^{(i,\gamma)} C_{\sigma_k \gamma, \sigma_k\dots\sigma_{0}}\cdot D^2_{\sigma_k\setminus i}\times S^1_{\gamma\cup i}=0.$$
On the other side,
\begin{multline*}(\delta'\Gamma^{k+1}_{p,q})_{\sigma_k\setminus
i,\sigma_k
\dots\sigma_{m+1}\sigma_{m-1}\dots\sigma_0}=\\=(-1)^{p+q-k+1}\sum_{\substack{\sigma_{m+1}\subset\sigma_m\subset\sigma_{m-1}\\|\sigma_m|=q-m}}\sum_{|\gamma|=p-q+k}(-1)^{(i,\gamma)}
C_{\sigma_k \gamma, \sigma_k\dots\sigma_{0}}\cdot
D^2_{\sigma_k\setminus i}\times S^1_{\gamma\cup i},\end{multline*}
hence $(\delta'\Gamma^{k+1}_{p,q})_{\sigma_k\setminus i,\sigma_k
\dots\sigma_{m+1}\sigma_{m-1}\dots\sigma_0}=0.$ We have shown that
$\partial \Gamma^{k}_{p,q}=-\delta' \Gamma^{k+1}_{p,q}.$\hfill
$\Box$

It follows from Theorem \ref{th.ret} and the construction of the
cell decomposition of the moment-angle complex
$\mathcal{Z}_\mathcal{K}$ that any cycle $\Gamma_s\in
H_s(\mathbb{C}^n\setminus Z_{\mathcal{K}})$ can be represented as a
sum of the cycles $\Gamma_{p,q}:$
$$\Gamma^s=\sum_{\substack{p+q=s\\p\geq q}}\Gamma_{p,q},$$
where $\Gamma_{p,q}$ is the cycles of the form
(\ref{eq.cycle_type}). From Lemma \ref{lemma.resolv} we have the
construction of the resolvent $\Gamma^0_{p,q},\dots,\Gamma^q_{p,q}$
of the cycle $\Gamma_{p,q}.$

The cover $\mathcal{U}_\mathcal{K}$ is
$\overline{\partial}$-acyclic. Indeed, all elements of the cover and
their intersections are isomorphic to
$\mathbb{C}^{n-k}\times(\mathbb{C^*})^{k}$ for appropriate choice of
$k$ and consequently are a Stein manifolds.

From Theorem \ref{th.BT1} and Theorem \ref{th.BT2} we obtain that
$H^s(K^\bullet(\mathcal{U}_\mathcal{K},\Omega^\bullet),D)$ is
isomorphic to the de Rham cohomology group
$H^s(\mathbb{C}^n\setminus Z_\mathcal{K},\mathbb{C}).$ Recall that
we use the following notation for the inclusions of complexes
$$\varepsilon:\mathcal{E}^\bullet(\mathbb{C}^n\setminus Z_\mathcal{K}) \rightarrow K^\bullet(\mathcal{U}_\mathcal{K},\mathcal{E}^\bullet),$$
$$\tau:K^\bullet(\mathcal{U}_\mathcal{K},\Omega^\bullet)\rightarrow K^\bullet(\mathcal{U}_\mathcal{K},\mathcal{E}^\bullet).$$
We will use the same notation for the induced isomorphisms on the
cohomology groups:
$$H^s(\mathbb{C}^n\setminus Z_\mathcal{K},\mathbb{C})\stackrel{\varepsilon}{\simeq} H^s(K^\bullet(\mathcal{U}_\mathcal{K},\mathcal{E}^\bullet),D)\stackrel{\tau}{\simeq} H^s(K^\bullet(\mathcal{U}_\mathcal{K},\Omega^\bullet),D).$$

\begin{lemma}\label{lemma.cocycle} Let $\omega\in  H^s(K^\bullet(\mathcal{U}_\mathcal{K},\Omega^\bullet),D),$ then there is a cocycle $\widetilde{\omega}$ such that $\widetilde{\omega}=\omega$ in $H^s(K^\bullet(\mathcal{U}_\mathcal{K},\Omega^\bullet),D),$
$\widetilde{\omega}=\sum_{p+q=s}\widetilde{\omega}^{p,q},$
$\widetilde{\omega}^{p,q}\in C^q(\mathcal{U}_\mathcal{K},\Omega^p),$
and
$$(\widetilde{\omega}^{p,q})_{\sigma_0,\dots,\sigma_q}=\sum_{|I|=p} C_{I,\sigma_0,\dots,\sigma_q} \frac{d z_I}{z_I},$$
 where $D \widetilde{\omega}^{p,q} =0$ for any $p$ and $q.$
\end{lemma}
{\it Proof.} Consider an arbitrary element $\omega$ of
$H^s(K^\bullet(\mathcal{U}_\mathcal{K},\Omega^\bullet),D),$ this
element is representable by cocycle
$\omega=\sum_{p+q=s}\omega^{p,q},$ where $\omega^{p,q}\in
C^q(\mathcal{U}_\mathcal{K},\Omega^p)$ and $\delta\omega^{p,q} = -d
\omega^{p-1,q+1}.$ The cocycle $\omega$ has a unique decomposition
$\omega= \widetilde{\omega} + \psi$, where
$(\widetilde{\omega}^{p,q})_{\sigma_0,\dots,\sigma_q}$ is the
following form
\begin{equation}\label{eq.cocycl}(\widetilde{\omega}^{p,q})_{\sigma_0,\dots,\sigma_q}=\sum_{|I|=p} C_{I,\sigma_0,\dots,\sigma_q} \frac{d z_I}{z_I},\end{equation}
and the Laurent expansion of
$(\psi^{p,q})_{\sigma_0,\dots,\sigma_q}$ does not contain summands
$\frac{d z_I}{z_I}.$

Let us show that $\widetilde{\omega}$ and $\psi$ are cocycles. Since
$\omega$ is a cocycle, we have
$\delta\widetilde{\omega}^{p,q}+\delta \psi^{p,q} = -d
\widetilde{\omega}^{p-1,q+1}-d \psi^{p-1,q+1}$. The forms $\frac{d
z_I}{z_I}$ are closed, hence $d \widetilde{\omega}^{p-1,q+1}=0.$
Since the Laurent expansions of the components of the cochain
$\delta \psi^{p,q}$ do not contain summand $\frac{d z_I}{z_I}$ and
the  cochain $\delta\widetilde{\omega}^{p,q}$ can be exact if and
only if $\delta\widetilde{\omega}^{p,q}=0$ (because nonzero linear
combinations of the forms $\frac{d z_I}{z_I}$ are nonexact on any
elements of the cover $\mathcal{U}_\mathcal{K}$),
$\delta\widetilde{\omega}^{p,q}=0.$ We get that
$\delta\widetilde{\omega}^{p,q}=d \widetilde{\omega}^{p,q}=0,$
consequently, $\widetilde{\omega}$ is a cocycle. The cochain
$\psi=\omega-\widetilde{\omega}$ is a difference of two cocycles,
hence $\psi$ is a cocycle.

Now we going to show that $\psi$ is a coboundary.

\begin{lemma}\label{lemma.zero}Let $\Gamma\in H_s(\mathbb{C}\setminus Z_\mathcal{K}),$ then $\int_{\Gamma}\varepsilon^{-1}\circ\tau(\psi)=0.$
\end{lemma}
{\it Proof.} For the cycle $\Gamma$ we have the expansion to the sum
$\Gamma=\sum_{\substack{p+q=s}}\Gamma_{p,q}.$ Lemma
\ref{lemma.resolv} gives the explicit construction of the resolvent
$\Gamma_{p,q}^0,\dots,\Gamma_{p,q}^q$ of the cycle $\Gamma_{p,q}.$
Let us construct a cocycle $\widetilde{\psi}_k=\sum_{p+q=s}
\widetilde{\psi}_k^{p,q}$ cohomologous to $\tau(\psi)$ of
$K^\bullet(\mathcal{U}_\mathcal{K},\mathcal{E}^\bullet)$

\begin{equation*}
\widetilde{\psi}_k^{p,q}=\left\{
           \begin{array}{ll}
             \psi^{p,q} & \mbox{\emph{for}\;} $\emph{q<k,}$ \\
             \psi^{p,q}-d \delta^{-1}(\psi^{p-1,q+1}-d \delta^{-1}(\psi^{p-2,q+2}-d \delta^{-1}(\dots-d \delta^{-1} (\psi^{0,p+q})))) & \mbox{\emph{for}\;} $\emph{q=k,}$\\
             0 & \mbox{\emph{for}\;} $\emph{q>k.}$ \\
           \end{array}
         \right. \end{equation*}

From Proposition \ref{prop.pairing} we obtain
$$\int_{\Gamma}\varepsilon^{-1}\circ\tau(\psi)=\sum_{p+q=s}\sum_{k\leq q}\langle\Gamma^k_{p,q},\widetilde{\psi}_q^{s-k,k}\rangle.$$
Let $k<q$ and $\omega\in\Omega^{p+q-k}(\mathcal{U}_{\sigma'}),$ it
is easy to see, that $\omega|_{D^2_{\sigma}\times S^1_{\gamma}}=0$
for \linebreak $|\sigma|=q-k>0,$ $|\gamma|=p-q+k,|$ $\sigma\subseteq
\sigma'.$ The forms
$(\widetilde{\psi}_q^{s-k,k})_{\alpha_0,\dots,\alpha_k}$ are
holomorphic on
$\mathcal{U}_{\alpha_0}\cap\dots\cap\mathcal{U}_{\alpha_k},$ indeed,
$(\widetilde{\psi}_q^{s-k,k})_{\alpha_0,\dots,\alpha_k}=(\psi^{s-k,k})_{\alpha_0,\dots,\alpha_k},$
on the other side, $(\Gamma^k_{p,q})_{\alpha_0,\dots,\alpha_k}$ is a
linear combination of the chains $D^2_{\sigma}\times S^1_{\gamma},$
$|\sigma|=q-k>0,|\gamma|=p-q+k$, combining these two facts we get
$\langle\Gamma^k_{p,q},\widetilde{\psi}_q^{s-k,k}\rangle=0.$

Consider the case $k=q.$ From the definition of $\widetilde{\psi}$
it follows that $\widetilde{\psi}_q^{p,q}=\psi^{p,q}+d \varphi$ for
some $\varphi\in C^{q}(\mathcal{E}^{p-1},\mathcal{U}).$  Since
$(\Gamma^q_{p,q})_{\alpha_0,\dots,\alpha_k}$ is a linear combination
of the cycles $S^1_{\gamma},|\gamma|=p,$
$$S^1_\gamma=\{|z_j|=1: j\in\gamma,z_k=1: k\not\in\gamma\},$$
 $\langle\Gamma^q_{p,q},\widetilde{\psi}_q^{s-q,q}\rangle=\langle\Gamma^q_{p,q},\psi^{s-q,q}\rangle.$ Indeed, by the Stokes formula  $\int_{S^1_{\gamma}} d \varphi=0,$ hence $\langle\Gamma^q_{p,q},d \varphi\rangle=0.$

Expand the from $(\psi^{s-q,q})_{\alpha_0,\dots,\alpha_q}$ to the
Laurent series,
$$(\psi^{s-q,q})_{\alpha_0,\dots,\alpha_q}=\sum_{a=(a_1,\dots,a_n)\in \mathbb{Z}^n} \sum_{|I|=p} C_{a ,I, \alpha_0\dots\alpha_q} z_1^{a_1}\dots z_n^{a_n} \frac{d z_I}{z_I}.$$
The integral $$\int_{S^1_{\gamma}} z_1^{a_1}\dots z_n^{a_n} \frac{d
z_I}{z_I}$$ is nonzero only if $I=\gamma$ and $a=0,$ i.e., for the
forms $\frac{d z_\gamma}{z_\gamma},$ but by the construction of
$\psi^{s-q,q}$ the Laurent expansion of
$(\psi^{s-q,q})_{\alpha_0,\dots,\alpha_q}$ does not contain the
summands   $\frac{d z_\gamma}{z_\gamma}.$ Consequently,
$\langle\Gamma^q_{p,q},\psi^{s-q,q}\rangle=0.$

We have shown that $\int_{\Gamma}\varepsilon^{-1}\circ\tau(\psi)=0.$
Lemma \ref{lemma.zero} is proved. \hfill $\Box$

By the de Rham Theorem any closed form $\omega$ of degree $s$ on
$\mathbb{C}^n\setminus Z_\mathcal{K}$ is exact if and only if
$\int_\Gamma \omega=0$ for any cycle $\Gamma\in
H_s(\mathbb{C}^n\setminus Z_\mathcal{K}).$ It follows from Lemma
\ref{lemma.zero} that $\varepsilon^{-1}\circ\tau(\psi)$ is
cohomologous to zero, hence $\psi$ is a coboundary . Lemma
\ref{lemma.cocycle} is proved.\hfill $\Box$

By Lemma \ref{lemma.cocycle} any cocycle $\omega\in
H^s(K^\bullet(\mathcal{U}_\mathcal{K},\Omega^\bullet),D)$ is
cohomologous to \linebreak
$\widetilde{\omega}=\sum_{p+q=s}\widetilde{\omega}^{p,q},$ where
$\widetilde{\omega}^{p,q}$ is of the form (\ref{eq.cocycl}).
Moreover, $\widetilde{\omega}^{p,q}\in
C^q(\mathcal{U}_\mathcal{K},\Omega^p)$ is a cocycle, i.e., $D
\widetilde{\omega}^{p,q}=0.$ We denote by
$H^{p,q}(K^\bullet(\mathcal{U}_\mathcal{K},\Omega^\bullet),D)$ a
subspace of
$H^s(K^\bullet(\mathcal{U}_\mathcal{K},\Omega^\bullet),D)$ generated
by cocycles $ \widetilde{\omega}^{p,q}.$  We obtain
$$H^s(K^\bullet(\mathcal{U}_\mathcal{K},\Omega^\bullet),D)=\bigoplus_{p+q=s} H^{p,q}(K^\bullet(\mathcal{U}_\mathcal{K},\Omega^\bullet),D).$$

Then the filtration $F^k
H^s(K^\bullet(\mathcal{U}_\mathcal{K},\Omega^\bullet),D)$ equals
$$F^k H^s(K^\bullet(\mathcal{U}_\mathcal{K},\Omega^\bullet),D)=\bigoplus_{p\geq k} H^{p,s-p}(K^\bullet(\mathcal{U}_\mathcal{K},\Omega^\bullet),D).$$
Hence, $$F^k H^s(\mathbb{C}^n\setminus
Z_{\mathcal{K}},\mathbb{C})\stackrel{\varepsilon^{-1}\circ
\tau}{\simeq}\bigoplus_{p\geq k}
H^{p,s-p}(K^\bullet(\mathcal{U}_\mathcal{K},\Omega^\bullet),D).$$

By the same argument as in Lemma  \ref{lemma.zero}  we obtain that
for every cycle $\Gamma_{p,q}\in H_{p,q}(\mathbb{C}^n\setminus
Z_{\mathcal{K}}),$ and every cocycle
$\widetilde{\omega}^{p',q'}\in
H^{p',q'}(K^\bullet(\mathcal{U}_\mathcal{K},\Omega^\bullet),D),$ the
following equality holds
$$\int_{\Gamma_{p,q}}\varepsilon^{-1}\circ\tau(\widetilde{\omega}^{p',q'})=0,$$
for $p\neq p',$ $q\neq q'.$ It follows from nondegeneracy of the
pairing between cohomology and homology that the pairing between
elements of $H_{p,q}(\mathbb{C}^n\setminus
Z_{\mathcal{K}},\mathbb{C})$ and \linebreak
$\varepsilon^{-1}\circ\tau
(H^{p',q'}(K^\bullet(\mathcal{U}_\mathcal{K},\Omega^\bullet),D))$ is
nondegenerate if $p=p',$ $q=q'$ and equals to zero otherwise. Thus,
$\varepsilon^{-1}\circ\tau(H^{p',q'}(K^\bullet(\mathcal{U}_\mathcal{K},\Omega^\bullet),D))=\phi(H^{p',q'}(R_\mathcal{K}\otimes\mathbb{C})).$
\hfill $\Box$

\section{Integral representations of holomorphic functions}

In the last section we study integral representations of holomorphic
functions such that kernels of these integral representations have
singularities on coordinate subspace arrangements in $\mathbb{C}^n$.
The examples of such integral representations are the
multidimensional Cauchy integral representation, whose kernel has
singularity on $(\{z_1=0\}\cup\dots\cup\{z_n=0\}),$ and the
Bochner–Martinelli integral representation, whose kernel has
singularity on $\{0\}$. In \cite{Sh} a family of new integral
representations of this kind was obtained, the kernels of these
integral representations have singularities on the subspace
arrangements defined by simple polytopes.


Denote by $U$ the unit polydisc in $\mathbb{C}^n:$
$$U=\{z=(z_1,\dots,z_n)\in\mathbb{C}^n: |z_i|<1, i=1,\dots,n\}.$$
Notice that the moment-angle complex $\mathcal{Z}_\mathcal{K}$ is
lying on the boundary  $\partial U$ of the polydisc.

\begin{theorem}
Given a nontrivial element $\omega'$  from $F^n
H^s(\mathbb{C}^n\setminus Z_{\mathcal{K}},\mathbb{C}).$ Then there
exists a closed $(n,s-n)$-form $\omega,$ $[\omega]=\omega'$ and an
$s$-dimensional cycle $\Gamma$ in $\mathbb{C}^n\setminus
Z_{\mathcal{K}}$ with support in $\mathcal{Z}_\mathcal{K},$ such
that for any function $f$ holomorphic in some neighborhood of $U$
the following integral representation holds
$$f(\zeta)=\int_{\Gamma}f(z) \omega(z-\zeta)$$ for $\zeta\in U.$\end{theorem}

{\it Proof.}
Since $\omega'\in F^n
H^s(\mathbb{C}^n\setminus~Z_\mathcal{K},\mathbb{C}),$ by Theorem
\ref{th.main} there is a cycle $\Gamma\in H_s(\mathbb{C}^n\setminus
Z_\mathcal{K},\mathbb{C}),$
$$\Gamma=\sum_{\substack{|\sigma|=s-n\\|\gamma|=2n-s}}C_{\sigma \gamma}\cdot D^2_\sigma\times S^1_\gamma,$$ such that $\langle\Gamma, \omega'\rangle =1.$
It follows from Lemma \ref{lemma.cocycle} that there exists a
cocycle $\omega^{n,s-n}\in
C^{s-n}(\mathcal{U}_\mathcal{K},\Omega^n),$
$$(\omega^{n,s-n})_{\alpha_0,\dots,\alpha_{s-n}}=B_{\alpha_0,\dots,\alpha_{s-n}} \frac{d z_1}{z_1}\wedge\dots\wedge\frac{d z_n}{z_n},$$ that is cohomologous to $\tau^{-1}\circ\varepsilon (\omega')$ in $H^s(K^\bullet(\mathcal{U}_\mathcal{K},\Omega^\bullet),D).$ The form \linebreak $\omega=\varepsilon^{-1}\circ(-\overline{\partial}\delta^{-1})^{s-n} \omega^{n,s-n}$ is cohomologous to $\omega',$ so
$$\int_{\Gamma} \omega =1.$$

Let us show that $\omega$ and $\gamma$ define an integral
representation. Consider the integral
$$\int_{\Gamma} \omega(z-\zeta),$$
where $\zeta\in U,$ here the notation $\omega(z-\zeta)$ stands for
the form $\omega$ after the change of coordinates $z \rightarrow
z-\zeta.$ Notice that the form $\omega(z-\zeta)$ is closed in $U,$
thus the integral of this form depends only on the homological class
of the integration cycle. Let us make a change of coordinates
$$\int_{\Gamma{-\zeta}}\omega(z),$$
 where $\Gamma{-\zeta}$ is a cycle $\Gamma$ shifted by the vector $-\zeta.$ In the sequel we will use the subindex $-\zeta$ to denote chains, cycles, and sets in  $\mathbb{C}^n$ shifted by the vector $-\zeta.$

Let us show that $\Gamma-\zeta$ is homologous to $\Gamma.$ Notice
that $(\mathcal{Z}_\mathcal{K}-\zeta)\cap Z_\mathcal{K}=\emptyset$
for any $\zeta\in U.$ Indeed,
$$\mathcal{Z}_\mathcal{K}-\zeta=\bigcup_{\sigma\in \mathcal{K}} (D^2_\sigma\times S^1_{[n]\setminus \sigma}-\zeta),$$
$$Z_\mathcal{K}=\bigcup_{\sigma\not\in\mathcal{K}}L_\sigma,$$
we see that $(D^2_\sigma\times S^1_{[n]\setminus \sigma}-\zeta)\cap
L_{\sigma'}=\emptyset$ for any $\sigma\in K, \sigma'\not\in K$ and
$\zeta\in U.$ Consider the chain $$\widetilde{\Gamma_{-\zeta}}=\{y:
y=x - t \zeta, x\in\Gamma,t\in[0,1]\},$$ the support of the chain
$\widetilde{\Gamma_{-\zeta}}$ is a subset of $\bigcup_{t\in
[0,1]}(\mathcal{Z}_\mathcal{K}-t\zeta),$  therefore
$\widetilde{\Gamma_{-\zeta}}$ is a subset of $\mathbb{C}^n\setminus
Z_\mathcal{K}.$
  Its boundary equals $\partial\widetilde{\Gamma_{-\zeta}}=(\Gamma-\zeta)-\Gamma,$ i.e., $(\Gamma-\zeta)$ and $\Gamma$ are homologous. So we have returned to the case $\int_{\Gamma}\omega(z),$ which was already considered.
We get
$$\int_{\Gamma}\omega(z-\zeta)=\int_{\Gamma{-\zeta}}\omega(z)=1.$$

By dentition $\omega(z-\zeta)$ is an $(n,s-n)$-form. Let $f(z)$ be a
function holomorphic in some neighborhood of unit polydisc $U$.
Since the operators $\overline{\partial}$ and $\delta$ are
interchangeable with the multiplication by a holomorphic function,
we get $f(z)
\cdot\omega(z-\zeta)=\varepsilon^{-1}\circ(-\overline{\partial}\delta^{-1})^{s-n}
f(z) \cdot\omega^{n,s-n}(z-\zeta).$ By Lemma \ref{lemma.resolv}
there is a resolvent $\Gamma^0,\dots,\Gamma^{s-n}$ of the cycle
$\Gamma$ such that
$$\Gamma^{s-n}_{\alpha_0,\dots,\alpha_{s-n}}=C'_{\alpha_0,\dots,\alpha_{s-n}}\cdot S^1_{[n]},$$
  $$S^1_{[n]}=\{|z_1|=\dots=|z_n|=1\}.$$
  Since $$\int_\Gamma \omega(z-\zeta)= \langle\Gamma^{s-n},\omega^{n,s-n}(z-\zeta)\rangle=1,$$ from the Cauchy integral representation formula we get $$\int_\Gamma f(z) \omega(z-\zeta)=\langle\Gamma^{s-n},f(z)\cdot\omega^{n,s-n}(z-\zeta)\rangle=f(\zeta).$$
\hfill $\Box$

\end{document}